\documentclass[11pt]{amsart} 
\usepackage{latexsym} 
\usepackage{amsfonts} 

\usepackage{amsmath,amsthm,amssymb}

\title[Quantitative aspects]{Some Quantitative Aspects of Fractional  Computability}

\author[I.~Kapovich]{Ilya Kapovich} 
 
\address{\tt Department of Mathematics, University of Illinois at 
  Urbana-Champaign, 1409 West Green Street, Urbana, IL 61801, USA 
  \newline http://www.math.uiuc.edu/\~{}kapovich/} \email{\tt 
  kapovich@math.uiuc.edu}

\author[P.~Schupp]{Paul Schupp} 
 
\address{\tt Department of Mathematics, University of Illinois at 
  Urbana-Champaign, 1409 West Green Street, Urbana, IL 61801, USA 
  \newline http://www.math.uiuc.edu/People/schupp.html} 
\email{schupp@math.uiuc.edu}

\newtheorem{thm}{Theorem}[section] 
\newtheorem{lem}[thm]{Lemma} 
\newtheorem{cor}[thm]{Corollary}  
\newtheorem{prop}[thm]{Proposition} \theoremstyle{definition} 
\newtheorem{defn}[thm]{Definition} 
\newtheorem{notation}[thm]{Notation} 
\newtheorem{conv}[thm]{Convention}

\begin{document} 
 
\begin{abstract} 
  In this article we apply the ideas of effective Baire category and 
  effective measure theory to study complexity classes of functions 
  which are ``fractionally computable" by a partial algorithm.  For 
  this purpose it is crucial to specify an allowable effective 
  density, $\delta$, of convergence for a partial algorithm.  The set 
  $\mathcal{FC}(\delta)$ consists of all total functions $ f: 
  \Sigma^\ast \to \{ 0,1 \}$ where $\Sigma$ is a finite alphabet with 
  $|\Sigma| \ge 2$ which are ``fractionally computable at density 
  $\delta$".  The space $\mathcal{FC}(\delta) $ is effectively of the 
  second category while any fractional complexity class, defined using 
  $\delta$ and any computable bound $\beta$ with respect to an 
  abstract Blum complexity measure, is effectively meager.  A 
  remarkable result of Kautz and Miltersen shows that relative to an 
  algorithmically random oracle $A$, the relativized class 
  $\mathcal{NP}^A$ does not have effective polynomial measure zero in 
  $\mathcal{E}^A$, the relativization of strict exponential time.  We 
  define the class $\mathcal{UFP}^A$ of all languages which are 
  fractionally decidable in polynomial time at ``a uniform rate'' by 
  algorithms with an oracle for $A$.  We show that this class does 
  have effective polynomial measure zero in $\mathcal{E}^A$ for every 
  oracle $A$. Thus relaxing the requirement of polynomial time 
  decidability to hold only for a fraction of possible inputs does not 
  compensate for the power of nondeterminism in the case of random 
  oracles. 
\end{abstract} 
 
\thanks{Both authors were supported by the NSF grant DMS-0404991. The first author is also supported by the NSF grant DMS-0603921} 
 

\maketitle

\section{Introduction}\label{intro} 
 
We now know that ``worst-case'' complexity measures such as 
polynomial time do not necessarily give a good overall picture of a 
particular problem or algorithm since it depends on the difficulty 
of the hardest instances of the problem, and these may be very 
sparse.  The  famous classic example of this phenomenon is Dantzig's 
Simplex Algorithm for linear programming.  The examples of V.~Klee 
and G.~Minty \cite{K-M} showing that the simplex algorithm can be 
made to take exponential time are very special.  A ``generic" or 
``random" linear programming problem is not ``special", and 
Dantzig's algorithm works quickly. Indeed, later algorithms which 
are provably polynomial-time have not replaced the simplex algorithm 
in practice.

Observations of this type led to the development of {\it average-case 
  complexity} by Gurevich ~\cite{Gur} and Levin~\cite{Lev}.  There are 
now different approaches to the average-case complexity, but they all 
require computing the expected value of the running time of an 
algorithm with respect to some measure on the set of inputs.  It is 
often difficult to establish an average-case result since a basic 
difficulty of worst-case complexity is still present: one needs a 
total algorithm which solves the problem and some upper bound on its 
worst-case difficulty.

Kapovich, Myasnikov, Schupp and Shpilrain ~\cite{KMSS} introduced the 
notion of \emph{generic-case} complexity, which deals with the 
performance of an algorithm on ``most'' inputs and completely ignores 
what happens on the ``sparse'' set of other inputs. They applied the 
idea to the classic decision problems of group theory - the word and 
conjugacy problems - and found that the ``linear programming 
phenomenon'' is extremely widespread there.  An important aspect of 
generic-case complexity is that it allows us to work with the entire 
class of partial computable functions, which is the natural setting of 
the general theory of computability, and one can often prove 
generic-case complexity results about  problems where the worst 
case complexity is unknown. 
 
This paper grew out of our interest in generic-case complexity but 
here we are interested in the more general concept of ``fractionally 
computable at an allowable density $\delta$.''  The basic idea is 
essentially the same as for generic-case complexity. However, we do 
not demand that the fraction of possible inputs on which a partial 
algorithm succeeds approaches one, but only that the algorithm 
succeeds at the given density $\delta$.

Specific questions about fractional complexity are important in 
cryptography, where one needs the assumption that problems such as 
calculating the discrete logarithm are generically difficult. 
Proposition $6.3$ in the book by Talbot and Welsh~\cite{TW} states 
that if there is a polynomial time algorithm which solves the discrete 
logarithm problem for a subset $B_p \subseteq \mathbb{Z}_p^{\ast}$ 
where $|B_p| \ge \epsilon |\mathbb{Z}_p ^{\ast}|$ then there is a 
probabilistic algorithm that solves the discrete logarithm problem in 
general with expected running time polynomial in $k$ and $1/\epsilon$.

In group theory, subgroups of finite index provide natural examples 
algorithms with fractional complexity.  Suppose that $G$ is a 
finitely generated group and $N$ is a normal subgroup of finite 
index $j$. Let $\psi$ be the natural homomorphism from $G$ onto 
$G/N$.  Then the algorithm for the word problem of $G$ which simply 
consists of answering ``no'' on input $w$ if $\psi(w) \neq 1$ works 
on the fraction $1 - \frac{1}{j}$ of inputs.  Note that there is no 
assumption about the complexity of the word problem for $G$.  (The 
same result holds for subgroups $H$ which are not normal by using 
coset diagrams.)  There are now several suggestions for using 
problems about various groups for the purposes of cryptography. 
Fractional computability issues, such as those coming from subgroups 
of finite index, may pose difficulties for security.

The main results of this paper are (see Sections 2 and 3 below for 
precise definitions): 
 
\begin{thm}\label{A} 
  For every Blum complexity measure $ {\Phi}$, for every allowable 
  density $\delta$ and for every effective bound $\beta$, the 
  fractional complexity class ${\Phi}[\beta,\delta]$ is effectively 
  meager in the space $\mathcal{FC}(\delta)$. 
\end{thm}

\begin{thm}\label{B} 
  For every oracle $A$ the set $\mathcal{UFP}^A$ has effective 
  polynomial-time measure zero with respect to $A$ in $\mathcal{E}^A$. 
\end{thm} 
 
In Theorem~\ref{A}, $\delta(n)$ is an \emph{effective density} of 
convergence for a partial algorithm. Roughly speaking, $\delta(n)$ 
specifies the fraction of all inputs of length $n$ in which a 
partial function under consideration is required to be defined.  The 
set $\mathcal{FC}(\delta)$ consists all total functions $ f: 
\Sigma^\ast \to \{ 0,1 \}$ where $\Sigma$ is a finite alphabet with 
$|\Sigma| \ge 2$ which are ``fractionally computable at density 
$\delta$". The function $\beta(n)$ is an effectively computable 
resource bound for  some abstract Blum complexity measure (e.g. 
time). Informally, the class ${\Phi}[\beta,\delta]$ consists of all 
partial computable functions that can be computed on a 
-fraction of the inputs of length $n$ which is at least $\delta(n)$ 
with a resource bound $\beta(n)$.  The space $\mathcal{FC}(\delta) $ is effectively 
of the second category while 
 any fractional complexity class, defined 
using $\delta$ and  resource bound $\beta$, is 
effectively meager. 
 
In Theorem~\ref{B}, the space $\mathcal{E}^A$ consists of all total 
functions computable in strict exponential time with an oracle for 
$A$. The space $\mathcal{UFP}^A$ consists of those functions in 
$\mathcal{E}^A$ that are partially calculable by partial computable 
functions that are \emph{uniform} and are computable in polynomial 
time.  Here a partial function $\phi$ from $\Sigma^\ast$ to $\{0,1\}$ 
is \emph{uniform} if there exists a positive integer $k$ such that for 
every $w\in \Sigma^\ast$ with $|w|\ge k$ there is some $z$ with 
$|z|\le k\log |w|$ such that $\phi(wz)$ is defined. Thus being uniform 
can be viewed as a version of ``fractional computability''.

A remarkable result of Kautz and Miltersen~\cite{KM} shows that for an 
algorithmically random set $A \subseteq \Sigma^*$ the class 
$\mathcal{NP}^A$ does not have effective polynomial-time measure zero 
in $\mathcal{E}^A$. Thus Theorem~\ref{B} above shows that fractional 
polynomial-time computability does compensate for  the power 
of nondeterminism.

The main lines of our considerations  are directly taken from 
known results in the theory of effective category and measure.  The 
 contribution of this paper consists in showing that 
 such results apply to the  study of 
fractional complexity.  classes.  We are particularly indebted to the 
book ``Computational Complexity: A Quantitative View'' by Marius 
Zimand~\cite{Z} and the articles by Calude~\cite{Cal} and by Kautz and 
Miltersen~\cite{KM}.

\section{Allowable Densities and Fractional Computability} 
\label{se:1-2} 
 
\begin{conv}\label{compl} 
 
  We \emph{fix} a finite alphabet $\Sigma$ with $k \ge 2$ letters 
  together with a linear ordering of the letters.  As usual, 
  $\Sigma^*$ denotes the set of all words over $\Sigma$.  If $w \in 
  \Sigma^*$ then the \emph{length}, $|w|$, of $w$ is the number of 
  letters in $w$.  We denote the empty word by $\lambda$.  The 
  \emph{canonical} or \emph{shortlex} ordering of $\Sigma^*$ lists 
  words in order of increasing length and within a given length by the 
  lexicographical order induced by the given alphabetical ordering of 
  $\Sigma$.  Thus for $\Sigma = \{a,b\}$ the list is 
\[ \lambda,  a, b, aa, ab, ba, bb,  aaa, .... \] 
We take the listing 
\[ w_1,w_2, w_3, ...\] 
as defining a bijection between $\Sigma^*$ and the natural numbers 
$\mathbb{N} = \{ 1,2,3,.....\}$.  Using this bijection, we can 
consider functions from $\Sigma^*$ to $\{0,1\}$ as functions from 
$\mathbb{N}$ to $\{0,1\}$.  In this article $\mathcal{F}$ denotes the 
set of total functions from $\Sigma^*$ to $\{0,1\}$. 
 
A \emph{language} $L$ over $\Sigma$ is a subset of $\Sigma^*$.  We can 
identify a language $L \subseteq \Sigma^*$ with its characteristic 
function $\chi_L$ where 
 
\[ 
\chi_L(n)= 
\begin{cases} 
  1 &\text{ if  $ w_n \in L$},\\ 
  0 &\text{ if $ w_n \not\in L$}. 
\end{cases} 
\] 
 
This identification gives a bijection between the set $\mathcal{L}$ of 
all languages over $\Sigma$ and the set $\mathcal{F}$ and we take 
these sets as being essentially the same.

A function $f$ from 
\[ \Sigma^* = \{ w_1,w_2,...,w_n,...\}  \] 
to $\{ 0,1 \}$ is an infinite sequence $(b_n)$ of 0's and 1's. If $f 
\in \mathcal{F}$ takes the value $1$ infinitely often we can regard 
$f$ as the unique binary expansion of a real number in the half-open 
unit interval $(0,1]$ which is not all $0$'s from some point onwards. 
\end{conv} 
 
Suppose that we have a \emph{partial algorithm} $\Omega$ for a set $S 
\subseteq \Sigma^*$.  In particular, this means that $\Omega$ is 
correct: If $\Omega$ converges on an input $w$ then $\Omega$ gives the 
correct answer as to whether or not $w \in S$.  We again point out 
that we completely ignore the performance of $\Omega$ on words not in 
$S$ and the complexity classes we consider will generally contain 
functions $f$ which are not computable.  Indeed, note that a single 
partial algorithm $\Omega$ generically computes uncountably many 
different functions if  the set $D$ on which the partial algorithm 
converges is generic while its complement $\overline D$ is infinite. 
Let $f'$ be the \emph{partial} function defined by $\Omega$.  Then we 
can choose values on the set $\overline D$ in a totally arbitrary way 
to complete $f'$ to a total function $f$ which is generically computed 
by the given algorithm.


\begin{conv} We want to fix an effective enumeration $(M_i)$ of all Turing machines 
  with input alphabet $\Sigma$ and with a special output tape consisting of 
  a single square in which a machine can print either $0$ or $1$.  Let 
  $\phi_i$ be the partial function from $\Sigma^*$ to $\{0,1 \}$ which 
  is computed by $M_i$. We write $\phi_i (x)\downarrow$ if $\phi_i$ 
  produces a value on input $x$. 
\end{conv}

A major concern of ~\cite{KMSS} was the rate of convergence of a given 
generic-case algorithm.  It turns out that this is not an accident and 
that a general discussion of fractional complexity classes requires 
providing an effective density function which specifies a lower bound 
on how many values must be defined at a given stage. We now regard the 
functions in $\mathcal{F}$ as functions $f:\mathbb{N} \to \{0,1\}$. 
We need ``the acceptable density so far'' to be defined at each input 
$n$. 
 
\begin{defn}  An \emph{allowable density function }  is a computable function 
  $\delta : \mathbb{N} \to \mathbb{Q} \cap [0,1]$ such that 
\[ 
\liminf_{n\to\infty} \delta(n)>0. 
\] 
 
Given an allowable $\delta$, if $\phi_i$ is a p.c. function, we write 
$\Delta(\phi_i)$ if the condition 
\[| \{m: m \le n, \phi_i(w_m) \downarrow \}|/n  \ge \delta(n) \] 
holds for all $n \ge 1$. 
\end{defn}

Note that there is no claim that the predicate $\Delta$ is computable. 
We can now precisely define the space $\mathcal{FC}(\delta)$ of 
functions which are fractionally computable at density $\delta$. We 
assume that an allowable density function $\delta$ is now fixed. 
 
\begin{notation} 
  If $\phi$ is a partial function and $f \in \mathcal{F}$ is a total 
  function, we write $\phi \sqsubseteq f$ if $\phi (x) = f(x)$ at all 
  arguments for which $\phi$ is defined. 
\end{notation}

\begin{defn} 
  Let $\delta$ be an allowable density function.  We define 
  $\mathcal{FC}(\delta)$ to be the space of all functions $f \in 
  \mathcal{F}$ such that \newline there exists a partial computable 
  function $\phi_i$ such that $\Delta(\phi_i)$ and $\phi_i \sqsubseteq 
  f$. 
\end{defn}

In order to define the appropriate topology on $\mathcal{FC}(\delta)$ 
we consider finite sequences $\tau = (s_1,\dots, s_n)$ where each $s_i 
$ is from the three-letter alphabet $ \{ 0, 1, \perp \}$. The symbol 
$\perp$ represents an undefined value.  If $\tau$ has length $n$ as a 
sequence we write $|\tau| = n$.  We also write $\tau(j)$ for $s_j$.

The set of positions for which $\tau$ is defined is 
 
\[ def(\tau)   = \{j:  j \le |\tau|, \tau(j) \ne \perp \}. \]

As for partial functions we write $\Delta(\tau)$ if $|\{j: j \le l, j 
\in def(\tau)|/l \ge \delta(l)$ for all $l \le |\tau|$. 
 
\begin{defn} 
  A finite sequence $\tau$ is $\delta$-\emph{allowable} if $\tau$ 
  contains at least one defined entry and $\Delta(\tau) $.  Let 
  $\mathcal{T}$ denote the set of all $\delta$-allowable finite 
  sequences. If $\tau_1$ and $\tau_2$ are allowable sequences we write 
  $\tau_1 \sqsubseteq \tau_2$ if $\tau_2$ agrees with $\tau_1$ at all 
  positions for which $\tau_1$ is defined. If $\tau \in \mathcal{T}$ 
  and $f \in \mathcal{F}$ we write $ \tau \sqsubseteq f $ if $f$ 
  agrees with $\tau$ at all positions at which $\tau$ is defined.

  If $\tau \in \mathcal{T}$ is an allowable sequence then the \emph{ 
    basic neighborhood defined by } $\tau$ is 
 \[  N(\tau) = \{ f :  f \in \mathcal{FC}(\delta) ,  \tau \sqsubseteq  f \}. \] 
 Note that if $\tau_1 \sqsubseteq \tau_2$ then $N(\tau_2) \subseteq 
 N(\tau_1)$ since $\tau_2$ specifies more information than $\tau_1$. 
 
\end{defn} 
 
Regarding sequences as words over the three-letter alphabet $\{ \perp, 
0,1 \}$ we can effectively enumerate all $\delta$-allowable sequences 
as 
\[ \tau_1, \tau_2, ..., \tau_n,...  \] 
by considering all finite sequences over $\{ \perp, 0,1 \}$ in the 
canonical order and successively listing only those sequences which 
are $\delta$-allowable.  It is easy to show that the collection $\{ 
N(\tau): \tau \in \mathcal{T} \}$ is a system of basic neighborhoods 
and we use the topology generated by this system.

\begin{prop}\label{prop:basic} 
  For every $\tau_i, \tau_j \in \mathcal{T}$ with $N(\tau_i) \cap 
  N(\tau_j) \ne \emptyset$ there exists $\tau \in \mathcal{T}$ with 
  $N(\tau) \subseteq N(\tau_i) \cap N(\tau_j) $. 
\end{prop} 
 
\begin{proof}   Since  $N(\tau_i) \cap  N(\tau_j) \ne \emptyset$, it follows that the sequences 
  $\tau_i$ and $ \tau_j$ agree at all positions where both are 
  defined.  Thus the following sequence $\tau$ of length $r = max\{ 
  |\tau_i|, |\tau_j| \}$ is well-defined.  For $x \le r$ let

\[ \tau(x) = 
\begin{cases} 
  &    \tau_i(x) \text{ if } x \in def(\tau_i), \\ 
  &    \tau_j(x)  \text{ if } x \in def(\tau_j), \\ 
  & \perp \ \text{ if} \ \tau_i (x) = \tau_j (x) = \perp. 
\end{cases} 
\] 
 
Since $\tau$ is defined where either of the $\delta$-allowable 
sequences $\tau_i$ or $\tau_j$ are defined, the sequence $\tau$ is 
$\delta$-allowable and $N(\tau) \subseteq N(\tau_i) \cap N(\tau_j)$ by 
definition. 
\end{proof} 
 
       Blum\cite{Blu} gave a very general definition 
of an abstract complexity measure and we work in that context since 
the specific nature of the complexity measure is not important.

\begin{defn}  A \emph{Blum Complexity Measure}  of partially computable functions 
  is a partially computable function $ {\Phi}(i,x)$ satisfying the 
  following two axioms: 
 
\begin{enumerate} 
\item $ {\Phi}(i,x) \downarrow \iff \phi_i (x) \downarrow. $ 
 
\item The cost predicate 
 
\[ Cost(i,x,y)  = 
\begin{cases} 
  1  & \text{ \ if \ }   {\Phi}(i,x) \le y,  \\ 
  0 & \text{ \ otherwise}. 
\end{cases} 
\] 
 
is computable. 
\end{enumerate} 
\end{defn} 
 
The standard measures of deterministic time or space are certainly 
Blum complexity measures.  For the remainder of this section we assume 
that some Blum Complexity Measure $ \Phi $ is fixed. 
 
We can now define fractional complexity classes using the complexity 
measure $\Phi$ and the density $\delta$.  Recall that $\{ \phi_i \}_i$ 
is an effective enumeration of partial computable functions from 
$\Sigma^*$ to $\{ 0,1 \}$ and that we think of such functions as being 
given by Turing machines which can print only the symbols $0$ and $1$ 
on their special output tape.  We now need to consider functions from 
$\Sigma^*$ to $\{ \perp, 0,1 \}$ and think that such functions are 
given by Turing machines which can print $0, 1$ or $\perp$ on their 
output tape.

\begin{defn} Let $\beta$ be any total computable function, which we 
  will refer to as the \emph{effective bound}. 
 
The function $\phi_i$ \emph{strictly bounded by} $ \beta$, which we denote 
by $\phi_i ^{[ \beta]}$, is defined as follows.  We take the Turing 
machine $M$ for $\phi_i$ and obtain the Turing machine $M'$ by adding 
an initial subroutine which, on input $x$, calculates 
$Cost(i,x,\beta(x))$.  If this value is $0$, then either $\Phi(i,x)$ 
is undefined (and hence $\phi_i(x)$ is undefined) or $\Phi(i,x)$ is 
defined and the complexity $\Phi(i,x)$ on input $x$ exceeds 
$\beta(x)$. In either case,  if $Cost(i,x,\beta(x))=0$, $M'$ prints the 
value $\perp$. If $Cost(i,x,\beta(x))=1$ (so that both $\phi_i(x)$ and 
$\Phi(i,x)$ are defined and, in addition, $\Phi(i,x)$ is bounded by 
$\beta(x)$), $M'$ prints the value calculated by $M$ on input $x$. 
\end{defn}

This construction gives us an effective enumeration  of all the 
functions $\phi_{i} ^{[\beta]}$.  Note that if these are considered as 
functions from $\Sigma^*$ to $\{0,1,\perp \}$ then they are total 
computable functions. Finally we have 
 
\begin{defn} 
  The fractional complexity class, $\Phi [ \beta,\delta ] $, defined 
  by $ \Phi, \beta$ and $\delta$ is 
 
\[ 
\Phi [ \beta, \delta ] = \{ f \in \mathcal{FC}(\delta ) : \exists i 
[\Delta (\phi_i ^{[\beta]}) \ \text{ and } \ \phi_i^{[\beta]} 
\sqsubseteq f] \} 
\]

\end{defn}

We now turn to the notion of effective Baire category.  The 
requirement for a set $S \subseteq \mathcal{FC}(\delta)$ to be 
effectively nowhere dense is that there is a uniform effective method 
which, when given any basic open neighborhood $N$, produces another 
basic neighborhood $N' \subseteq N$ such that $S \cap N' = \emptyset$. 
For a meager set, that is, a countable union of nowhere dense sets, we 
require that the method be uniform over all the members of the union. 
Recall that $\mathcal T$ denotes the set of all $\delta$-allowable 
finite sequences. 
 
\begin{defn} 
 
  A set $X \subseteq \mathcal{FC}(\delta)$ is \emph{ effectively 
    nowhere dense } in $\mathcal{FC}(\delta)$ if there exists a total 
  computable witness function $\alpha: \mathcal{T} \to \mathcal{T} $ 
  such that: 
 
\begin{enumerate} 
\item $ \tau \sqsubseteq \alpha(\tau) \text{ for all } \tau \in 
  \mathcal T.$ 
 
\item $ X \cap N(\alpha (\tau))= \emptyset.$ 
\end{enumerate} 
 
A set $ X \subseteq \mathcal{FC}(\delta)$ is \emph{ effectively meager 
} if there exist a sequence of nowhere dense sets $(X_i)_i$ and a 
total computable witness function $\alpha: \mathbb{N} \times 
\mathcal{T} \to \mathcal{T} $ of two variables such that: 
\begin{enumerate} 
 
\item $ X = \bigcup_{i = 1}^\infty X_i$ 
 
\item $ \tau \sqsubseteq \alpha (i, \tau ) \text{ for all } (i,\tau) 
  \in \mathbb{N} \times \mathcal{T}.$ 
 
\item $ X_i \cap N(\alpha(i, \tau)) = \emptyset.$ 
\end{enumerate}

A set is \emph{effectively ample} (effectively of the second category) 
if it is not effectively meager. 
\end{defn} 
 
It is now easy to prove the desired result that any complexity class 
$\Phi [\beta, \delta]$ is effectively meager while the entire space 
$\mathcal{FC}(\delta)$ is effectively of the second category.  Indeed, 
we have the following result. (Compare ~\cite{Cal}.) 
 
\begin{lem} 
  For every meager set $X$ and for every $\tau \in \mathcal{T}$, there 
  is a total computable function $f \in N(\tau) - X$. 
\end{lem} 
 
\begin{proof} 
  Since $X$ is effectively meager we can write $X = \bigcup_{i = 
    1}^\infty X_i$ where $X$ is effectively meager via the witness 
  function $\alpha(i,\tau)$. We define a total computable function $f$ 
  iteratively by a simple diagonalization argument.  For a given 
  $\tau$ let $\sigma_0$ be the sequence of length $|\tau|+1$ agreeing 
  with $\tau$ at all places where $\tau$ is defined and having $0$ in 
  all places where $\tau$ is undefined, and with $0$ as the last entry 
  of the sequence $\sigma_0$.  Then $\tau \sqsubseteq \sigma_0$. 
 
  Let $\eta_1 = \alpha(1,\sigma_0)$.  Let $\sigma_1$ be the sequence 
  of length $|\eta_1| + 1$ which agrees with $\eta_1$ in all places 
  where $\eta_1$ is defined and which has $0$ in all places where 
  $\eta_1$ is undefined, and with $0$ as the last entry of the 
  sequence $\sigma_1$.  So $ | \sigma_1| > |\eta_1|$ and all entries 
  in $\sigma_1$ are defined.  Since $\eta_1 \sqsubseteq \sigma_1$ , 
  $N(\sigma_1) \cap X_1 = \emptyset$. 
 
  We continue in the same fashion. Let $\eta_2 = \alpha(2,\sigma_1)$. 
Let $\sigma_2$ be the sequence of length $|\eta_2| + 1$ which agrees  with $\eta_2$ in all places where 
$\eta_2$ is defined and which has  $0$ in all places where $\eta_2$ is undefined, and 
with  one more defined position with entry $0$ at the end of $\eta_2$.  Thus $|\sigma_2| = |\eta_2| + 1$ 
and all entries in $\sigma_2$ are defined. 
Since $\sigma_2 \sqsubseteq \eta_2$, we have $ N(\sigma_2) \cap X_2 = \emptyset$. 
 
By this process, we iteratively define a sequence $(\sigma_i)_i$ of 
$\delta$-allowable intervals $\sigma_i$ in which all entries are 
defined such that $N(\sigma_i) \cap X_i = \emptyset$ and such that 
$\sigma_i\sqsubseteq \sigma_{i+1}$ and $|\sigma_i|<|\sigma_{i+1}|$ for 
every $i$.  Let $\sigma=\sigma(1), \sigma(2),\dots$ be the infinite 
binary sequence such that for every $i$ the initial segment of 
$\sigma$ of length $|\sigma_i|$ is $\sigma_i$. Note that every initial 
segment of $\sigma$ is a $\delta$-allowable sequence and that 
$|\sigma_i|\ge i$ for every $i$. 
 
Consider the function $f$ defined as $f(j)=\sigma(j)$ for every $j$. 
Clearly, $f$ is a total computable function, since for every $n$ we 
have $n\le |\sigma_n|$ and $f(n)=\sigma_n(n)$. 
 
Now $f \in N(\tau)$, since it agrees with $\tau$ at all places where 
$\tau$ is defined, and $f \notin X_i$ for all $i$. 
 
\end{proof} 
 
The theorem immediately yields the following corollary. 
 
\begin{cor} The set $\mathcal{R}$ of total effectively computable  functions from  $\Sigma^*$ to $\{0,1\}$ 
  is not meager in the space $\mathcal{FC}(\delta)$. 
\end{cor} 
 
\begin{thm}\label{thm:A} 
  For every Blum complexity measure ${\Phi}$, for every allowable 
  density $\delta$ and for every effective bound $\beta$, the 
  fractional complexity class ${\Phi}[\beta,\delta]$ is effectively 
  meager in the space $\mathcal{FC}(\delta)$. 
\end{thm} 
 
\begin{proof} 
  We have an effective enumeration $\ (\phi_{i}^{[\beta]} )_i $ of all 
  strictly $\beta$-bounded partial functions.  Let 
\[ 
C_i = \begin{cases} 
  &  \{ f \in \mathcal{FC}(\delta): \phi_i^{[\beta]} \sqsubseteq f \} \text{ \ if  \ } \Delta(\phi_i^{[\beta]}) \\ 
  & \emptyset \ \text{otherwise}. 
         \end{cases} 
\]

It is clear that $\Phi[\beta,\delta] = \bigcup_{i} C_i $ so we need 
only specify an effective witness function $\alpha$. Given an index 
$i$ and a $\delta$-allowable sequence $\tau$ compute 
$\phi_i^{[\beta]}$ on the first $|\tau|$ inputs in the canonical 
order.  If the computed sequence $\sigma$ of length $|\tau|$ is not 
$\delta$-allowable then $C_i = \emptyset$ and we set $\alpha(i,\tau) = 
\tau$.  Suppose now that $\sigma$ is $\delta$-allowable.  If $\sigma$ 
has a defined value $v$ on an input $w_j$ with $v \ne \tau(j)$ again 
set $\alpha(i,\tau) = \tau$.  Suppose now that $\sigma$ is allowable 
and that for all $j\le |\tau|$ with a defined value $\sigma(j)$ we 
have $\sigma(j)=\tau(j)$. 
 
We claim that there exists $r>|\tau|$ such that either 
$\phi_i^{[\beta]}$ has a defined value $\phi_i^{[\beta]}(r)$ or the 
sequence $\sigma_r$ of the values of $\phi_i^{[\beta]}$ on the first 
$r$ inputs is non-allowable. This follows from the assumption 
$\liminf_{n\to\infty} \delta(n)>0$ in the definition of an allowable 
density function and from the definition of a $\delta$-allowable 
sequence. We continue computing values of  $\phi_i^{[\beta]}$ 
until we find  the smallest $r>|\tau|$ with the above property. 
 
If the sequence $\sigma_r$ is not allowable then $C_i=\emptyset$ and 
we again set $\alpha (i, \tau) = \tau$. If $\sigma_r$ is allowable and 
the $r$-th entry of $\sigma_r$ is a defined value $v$, we set $\alpha 
(i, \tau)$ to be the sequence agreeing with $\sigma_r$ at all the 
positions $j < r$ and having value $1-v$ at position $r$.  In either 
case we have $\alpha(i,\tau) \cap C_i = \emptyset$. 
\end{proof} 
 
Note that, in general, a fractional complexity class $ 
  \Phi[\delta,\beta]$ contains uncountably many functions while 
the nonmeager set $\mathcal{R}$ is countable.

\section{Nondeterminism versus fractional  polynomial-time computability} 
 
It should be expected that partial computability at a fixed density 
cannot make great inroads into the power of nondeterminism.  A 
nondeterministic machine can guess on \emph{every} input, while in 
considering fractional complexity, we still have a deterministic 
machine which is required to actually do the desired calculation on a 
non-negligible set of inputs. 
 
Turing himself~\cite{Tur} introduced the idea of Turing machines with 
an oracle.  We think of an oracle Turing machine as a Turing machine 
with a special hardware slot and any set $A \subseteq \Sigma^*$ can be 
``plugged into'' the slot.  The machine has a special query tape and a 
``branching instruction'' in addition to the standard Turing machine 
instructions.  The branching instruction has the form $ q_i , \sigma_l 
\rightarrow q_j , q_k$.   It is 
crucial that  all oracle machines are still 
specified by finite programs of instructions of the two types, 
so we still have an effective enumeration of all oracle 
Turing machines.  In a \emph{Turing machine} $M^A$ \emph{with an 
  oracle for} $A$, an instruction $q_i , \sigma_l \rightarrow q_j , q_k$ works as 
follows. If the machine $M^A$ is in state $q_i$ reading the symbol 
$\sigma_l$ on its work tape then the machine goes to state $q_j$ if 
the word written on the query tape belongs to the set $A$ and goes to 
state $q_k$ if the word on the query tape is not in the set $A$.

``Classical'' results of computability theory ``relativize'' in the 
following strong sense.  For example, take the proof of the 
unsolvability of the Halting Problem.  Not only the statement of the 
theorem but the given proof remain correct if one everywhere replaces 
the words ``Turing machine'' by the words ``Turing machine with an 
oracle for A''.  One could take this relativization property as a 
definition of ``classical''. 
 
However, the well-known theorem of Baker, Gill and Solovay,~\cite{BGS} 
showed that the question of $\mathcal{P}$ versus $\mathcal{NP}$ does 
\emph{not} relativize.  It is easy to construct an oracle $A$ such 
that $\mathcal{P}^A = \mathcal{NP}^A$.  Indeed, any set $A$ which is 
complete for $PSPACE$ will do.  But there are many oracles $A$ for 
which $\mathcal{P}^A \ne \mathcal{NP}^A$.Indeed, 
  Bennett and Gill~\cite{BG} showed that 
$\mathcal{P}^A \ne \mathcal{NP}^A$ with respect to a ``random'' 
oracle.  This means that  the set of $A$ such that $\mathcal{P}^A 
\ne \mathcal{NP}^A$ has Lebesgue measure one in the space of all 
languages over $\Sigma$.  Later results show that for a random oracle 
the separation between $\mathcal{P}^A$ and $ \mathcal{NP}^A$ is indeed 
very strong.  Our approach in this section is inspired by the 
remarkable result of Kautz and Miltersen~\cite{KM} which we will 
explain below. We use this approach to show that requiring polynomial time 
computation to succeed only on a ``reasonable fraction'' of the inputs 
does not significantly improve our computing power when compared to 
nondeterminism for ``algorithmically random'' oracles.

  First of all,  the ideas  of generic-case computability, 
and indeed fractional computability at an allowable density $\delta$ 
certainly  relativize without any problem.  All 
definitions are exactly the same except that we now 
consider Turing machines with an oracle for $A$.

For this section we work inside the class of functions 
 \[\mathcal{E}^A = \bigcup_{c} DTIME^A (2^{cn}+ c),\] 
 computable in 
strict exponential time by Turing machines with an oracle for $A$. 
Note that if we are working with respect to an oracle $A$ then the elements of 
$\mathcal{E}^A$ are total functions $\Sigma^\ast\to \{0,1\}$. 
 
Effective measure theory was formulated by Lutz~\cite{Lut} building on 
earlier work of Schnorr~\cite{Schn}.  Recall that in discussing 
$\mathcal{P}$ and $\mathcal{NP}$ we are considering sets of languages 
over an alphabet $\Sigma$.  As mentioned earlier, we identify a 
language $L$ with the infinite binary sequence specifying its 
characteristic function.  We have the canonical enumeration of all 
words $w_1,...,w_n,...$ of all words in $\Sigma^*$.  We think of $L$ 
as the infinite binary sequence $L(0),...,L(n),...$ where $L(n) = 1$ 
if $w_n \in L$ and $L(n) = 0$ otherwise.  We use the formulation of 
effective measure theory in terms of computable martingales, which are 
strategies for betting on the values of successive bits of an infinite 
binary sequence.  Formally, 
 
\begin{defn} 
  A \emph{martingale} is a function $d: \{0,1\}^* \to \mathbb{R}$ such 
  that for all $\sigma \in \{0,1\}^*$ 
\[   (*) \hspace{.1in}  d(\sigma) = \frac{d(\sigma 0) + d(\sigma 1)}{2}  \] 
and the value $d(\lambda)$ of $d$ on the empty word is greater than 
$0$.

The martingale $d$ \emph{succeeds} on a sequence $\alpha \in \{ 0,1 
\}^{\infty}$ if 
\[ \limsup_{n \to \infty}  d(\alpha[1,\dots ,n]) = \infty . \] 
where $\alpha[1,\dots,n]$ is the initial segment of $\alpha$ of length $n$. 
The martingale $d$ \emph{succeeds} on a set $S \subseteq \{ 0,1 
\}^{\infty}$ if it succeeds on all sequences in $S$. 
\end{defn} 
 
We can think that we start with one dollar and double the bet each 
time, splitting the bet between the two possible next values according 
to the strategy $d$ .  We succeed on the set $S$ if we win an infinite 
amount of money on every sequence in $S$.  If we think of 
$\{0,1\}^{\infty}$ as the unit interval one can show that a set $C 
\subseteq [0,1]$ has Lebesgue measure $0$ if and only there exists 
some martingale which succeeds on $C$. 
 
For effective measure theory one imposes a condition on the difficulty 
of computing a martingale. We are interested in martingales which are 
computable in polynomial time with respect to a fixed oracle $A$. 
 
\begin{defn} 
  An \emph{ $A$-polynomial-time martingale } is a function $d: 
  \{0,1\}^* \to \mathbb{Q}$ which satisfies the martingale equation 
  (*) and which is computable in polynomial time by by some Turing 
  machine with an oracle for the set $A$.

  A set $S \subseteq \{0,1\}^{\infty} $ has \emph{ effective 
    polynomial-time measure zero with respect to } $A$ if there exists 
  an $A$-polynomial-time martingale $d$ which succeeds on all 
  sequences in $S$.  We write ``$S$ has effective $P^A$ measure zero''. 
\end{defn} 
 
Recall that we are working inside a space $\mathcal{E}^A$ of functions 
computable in strict exponential time by Turing machines with an 
oracle for $A$. The argument given in Zimand~\cite{Z} relativizes to 
give: 
 
\begin{thm}\cite{Z}  The set $\mathcal{E}^A$ does not have effective $P^A$-measure zero. 
\end{thm} 
 
\begin{proof} 
  For every $A$-polynomial time martingale $d$ we define a language $L 
  \in \mathcal{E}^A$ on which $d$ does not succeed.  The martingale 
  equation $2 d(w) = d(w0) + d(w1) $ implies that either $d(w0) \le 
  d(w)$ or $d(w1) \le d(w)$.  We put the empty word $\lambda$ in $L$ 
  and then iteratively define $L$.  If $\sigma = L[1,..n]$ has already 
  been defined, then $L[1,...,n+1] = \sigma 1$ if $d(\sigma 1) \le 
  d(\sigma)$ and $L[1,...,n+1] = \sigma 0$ otherwise.  It is clear 
  that $d$ does not succeed on $L$ since $d(L[1,...,n]) \le 
  d(\lambda)$ for all $n$.

  We need only check that $L \in \mathcal{E}^A$. Given an arbitrary $w 
  \in \{ 0,1 \}^*$, with $|w| = n$, we possibly need to calculate $d$ 
  on all words of length of length $(n-1)$.  There is a constant $c$ 
  such that on inputs of length $r$ $d$ is calculable in time $r^c + 
  c$ by a Turing machine with an oracle for $A$. Thus the entire 
  calculation can be done in time $2^{n-1}[ (n-1)^c ]$ so $L \in 
  \mathcal{E}^A$.  $\square$ 
\end{proof}

In their remarkable article, Kautz and Miltersen~\cite{KM} use the 
concept of sets which are ``algorithmically random'' in the sense of 
Martin-Lof~\cite{ML}. The precise details of that definition need not 
to be given here and the important point for us is that it yields a 
large class of sets for which the following theorem of Kautz and 
Miltersen holds.

\begin{thm}[Kautz, Miltersen~\cite{KM}] 
  If $A \subseteq \Sigma^*$ is an algorithmically random set then the 
  set $\mathcal{NP}^A$ does not have effective $P^A$-measure zero in 
  $\mathcal{E}^A$. 
\end{thm}

In order to discuss fractional polynomial time computability we again 
need to impose a suitable effective density condition which now 
becomes ``uniformity''. 
 
\begin{defn} 
A partial function $\phi$ from $\Sigma^*$ to $\{ 0,1 \}$ is $k$-\emph{uniform} 
if for all $w \in \Sigma^*$ with $|w|  \ge k$,  there exists a $z$ with $|z| \le k \  log(|w|)$ such that 
$\phi(wz) \downarrow $.  Thus for every $w$ there is a ``reasonably short'' $z$ such 
that $\phi_i$ converges on $wz$. 
 
A partial function $\phi$ is \emph{uniform} if it is 
$k$-uniform for some positive integer $k$. We write $U( \phi)$ if 
$\phi$ is uniform. 
\end{defn} 
 
   Note  that if we have an algorithm $\Omega$ which 
generically solves a decision problem, then for every $w$ there is 
some $z$ such that $\Omega$ converges on $wz$.  This is because any 
cylinder $C = \{ wu \}$ consisting of all words with prefix $w$ is not 
a negligible set.

\begin{conv} 
From now on, we will assume that $\Sigma=\{0,1\}$ although all the 
arguments below work for an arbitrary finite alphabet $\Sigma$. 
 
In general, a superscript $A$ for a function, such as $\phi^A$, 
indicates that $\phi^A$ is a partial function computable by a Turing 
machine with an oracle for $A$. Similarly, a superscript $A$ for a 
Turing machine, such as $M^A$, indicates that $M^A$ is a Turing 
machine with an oracle for $A$. 
\end{conv}

 \begin{defn} 
   If $\phi_i^A$ is a partial computable function, computed by the 
   $i$-th Turing machine $M_i^A$ with an oracle for $A$, the function 
   $\phi_i^A [i]$ is the function computed as follows.  We modify $M_i^A$ 
   to a Turing machine $Q_i^A$ by adding a subroutine to force the 
  the function obtained to be $i$-uniform with its computation time 
bounded by $n^{ci}$  
   on inputs $|w|$ with $|w| \ge i$, where $c$ is a constant independent of 
   $i$ and $w$.

   In detail, on an input $w$, $Q_i^A$ prints $\perp$ if $|w| < i$. 
 
   Suppose now that $|w|\ge i$. Then $Q_i^A$ carries out the 
   computation of $M_i^A$ for $n^{i}$ steps.  If $M_i^A$ calculates a 
   value from $\{0,1\}$, then $Q_i^A$ prints that value.  If not, $Q_i^A$ considers, in 
   the canonical order, the extensions $wz$ with $|z| \le i \log 
   (|w|)$ and carries out the computation of $M_i^A$ on $wz$ for $|w|^i$ 
   steps.  If $M_i^A$ calculates a value on such an extension, then the 
   condition that we are calculating a $i$-uniform function is 
   verified for the input $w$ and $Q_i^A$ outputs the value $\perp$ for 
   input $|w|$.  If $M_i^A$ does not calculate a value on any of these 
   extensions, then $Q_i^A$ outputs the value $0$ for input $w$, again 
   ensuring that the calculated function is $i$-uniform. 
 
   The number of words $z$ with $|z| \le i \log (|w|)$ is 
   $2^i 2^{\log (|w|)}= 2^i |w|$. It follows that for every $w$ with 
   $|w|\ge i$ the machine $Q_i^A$ prints a value $0,1$ or $\perp$ in at 
   most $|w|^i+|w|^i2^i|w|$ steps. Recall, that if $|w|\le i-1$, then 
   $Q_i^A$ prints the value $\perp$ in the input $w$. Thus for every 
   $w\in \Sigma^\ast$ the machine $Q_i^A$ computes a value from 
   $\{0,1,\perp\}$ on the input $w$ in $\le |w|^{ci}$ steps where 
   $c>0$ is independent of $i$ and $w$.

   Since we uniformly effectively obtain $Q_i^A$ from $M_i^A$, there is an 
   effective enumeration of all the functions $\phi_i^A[i]$.  Note 
   that since any particular partial computable function has 
   infinitely many indices, for any partial computable function 
   $\phi^A$ which is $k$-uniform for some $k$ and whose computation 
   time on inputs for which it calculates a value is bounded by a 
   polynomial, there is a large enough index $i$ such that 
   $\phi_i^A[i](w) = \phi^A(w)$ for all inputs with $|w| \ge i$. 
\end{defn} 
 
Recall that we are identifying languages with their characteristic 
functions. 
 
\begin{defn} 
  We consider the set $\mathcal{UFP}^A$ of all those languages 
  (functions) in $\mathcal{E}^A$ which are partially calculable by 
  partial computable functions which are uniform with computation time 
  strictly bounded by a polynomial $n^j$ on some Turing machine with 
  an oracle for $A$. Formally, 
\[ 
 \mathcal{UFP}^A  = \{ f : f \in \mathcal{E}^A,  \text{ and there is 
   some $i$ such that } \phi_i^A [i] \sqsubseteq f  \} 
 \] 
 
\end{defn}

For the next theorem we essentially use the proof in section 3.4 of 
Zimand~\cite{Z} that polynomial time $\mathcal{P}$ has effective 
polynomial-time measure zero, noting that it applies to 
$\mathcal{UFP}^A$. 
 
\begin{thm}\label{thm:B} 
  For every oracle $A$, the set $\mathcal{UFP}^A$ has effective $P^A$- 
  measure zero in $\mathcal{E}^A$ . 
\end{thm}

\begin{proof} 
 
  First of all, as noted above, we can give an effective enumeration 
  $\{ {Q_i}^A \}$ of all Turing machines with an oracle for $A$ such 
  that ${Q_i}^A$ calculates $\phi_i^A [i]$.  This means that we have 
  one Turing machine $Q^A (i, w)$ such that for every $w\in 
  \Sigma^\ast$ $Q^A (i,w)$ simulates ${Q_i}^A$ on input $w$ in time 
  bounded by $(\log i)^{c_1} |w|^{c_2 i}$, where $c_1>0, c_2>0$ are 
  constants independent of $i$ and $|w|$.

  Let $S_i$ be the set of functions 
\[ S_i = \{ f \in \mathcal{E}^A:  \phi_i^A [i] \sqsubseteq f \}. \] 
Then $\mathcal{UFP}^A = \bigcup_{i \in \mathbb{N}} S_i$.  We define 
a martingale which succeeds on $\mathcal{UFP}^A$ in three stages. 
 
First, we need to define a martingale $d_i$ which succeeds on on the 
set $S_i$.  We use the variable $x$ to denote arguments to a 
martingale.  Since a martingale is betting on characteristic sequences 
of languages, the position $x(n)$ is supposed to tell us the value 
$f(w_n)$ for the functions $f$ in $S_i$.  It is important to 
keep in mind that $|w_n| \le \log n$. (By $ \log n$ we mean 
 $ \lfloor \log_2 n  \rfloor)$. ) 
 
Let $x\in \Sigma^\ast$ and let $n=|x|$. 
 
If $|w_n|<i-1$, we put $d_i(x)=1$.

Suppose that $|w_n|\ge i-1$. 
We set:

\[ 
   d_i (x 0):= 
    \begin{cases} 
      & \  2 d_i(x) \ \text{ if } \phi_i^A [i](w_{n+1}) = 0,\\ 
      & \  0 \ \text{ if } \phi_i^A[i](w_{n+1}) = 1,\\ 
      & \  d_i(x) \ \text{ if } \phi_i^A[i](w_{n+1}) = \perp,\\ 
     \end{cases} 
\] 
 
and 
 
\[ 
d_i (x 1):= 
\begin{cases} 
          & \  0 \ \text{ if } \phi_i^A [i](w_{n+1}) = 0,\\ 
          & \  2 d_i(x) \ \text{ if } \phi_i^A [i](w_{n+1}) = 1,\\ 
          & \  d_i(x) \ \text{ if } \phi_i^A [i](w_{n+1}) = \perp.\\ 
\end{cases} 
\] 
It is easy to see that $d_i$ is a martingale.

Suppose that $j_1<\dots <j_m$ are indices such that $|w_{j_s}|\ge i$ and 
$\phi_i^A[i](w_{j_s})$ is defined for $s=1,\dots m$. Then for any $x\in 
\Sigma^\ast$ with $|x|=j_m$ such that 
$\phi_i^A[i][1,\dots,j_m]\sqsubseteq x$ 
we have 
\[ 
d_i(x)=d_i(x(1)\dots x(j_1)\dots x(j_2)\dots x(j_m))\ge 2^m. 
\] 
Hence $d_i$ succeeds on $S_i$. 
There is a Turing machine $D^A(i,x)$ with an oracle for $A$, which, given $i$ and 
$x\in \{0,1\}^\ast$ with $|x|=n$, computes $d_i(x)$ in time bounded by 
\[ 
 n (\log i)^{c_1}  (\log n)^{c_2 i}. 
\] 
If $|x|=n\ge i$, this time is at most 
\[ 
n (\log n)^{c_1}  (\log n)^{c_2 i}\le n (\log n)^{c_3i}\tag{$\dag$} 
\] 
where $c_3>0$ is independent of $i,n$. 
Similarly, if $n\ge |x|$ and $n\ge i$ then $d_i(x)$ is computed in 
time bounded by the estimate $(\dag)$.

Second, in order to obtain a global martingale which is calculable 
in polynomial time we need to exponentially inflate indices.  Let 
$\widetilde{S}_{2^{2^i}} = S_i$,  and let $\widetilde{S}_j = 
\emptyset$ if $j$ does not have the form $2^{2^i}$. 
 
Let  $\widetilde{d}_j$ be the constant martingale assigning $1$ to 
all inputs if $S_j = \emptyset$. Let $\widetilde{d}_j  = d_i$ if $j 
= 2^{2^i}$. In the this  case,  $\widetilde{d}_j(x)$ can be 
calculated in time \[ \le n(\log n)^{c_3 \log(\log j)} \text{ for 
$|x|=n \ge j$}.\]

We now need the inequality \[(\log n )^{\log \log j} \le n \text{ 
for }n \ge j\ge 2.\tag{$\ddag$}\] Since $\log$ is an increasing 
function, if we fix $n$, it suffices to prove the inequality for $j 
= n$. Taking logs of both sides of the inequality 
and setting $j =n$ we need 
\[ (\log (\log n))^2 \le \log n \] 
which holds for $ n \ge 2$. 
 
The inequalities $(\dag)$ and $(\ddag)$ imply that 
$\widetilde{d_j}(x)$ can be calculated for $|x|=n \ge j$ in time 
$n^{c_4}$, where $c_4$ is independent of $j,n$. The same is true if 
$n\ge |x|$ and $n\ge j$.

Third, we now need to define another martingale $\widehat{d_j}$ 
which dampens  $\widetilde{d}$.

If $x \in \{ 0,1 \}^*$ is nonempty, let $pref(x)$ denote the prefix of 
$x$ of length $|x| - 1$.  Let $\delta_j (x) $ be defined for nonempty 
$x$ by $\widetilde{d}_j(x) = \delta_j (x) \widetilde{d}_j (pref(x))$, 
provided $\widetilde{d}_j (pref(x))\ne 0$.  If $\widetilde{d}_j 
(pref(x))=0$ , we put $\delta_j (x)=1$. 
 
Note that if $\widetilde{d}_j (pref(x))=0$ then $\widetilde{d}_j 
(x)=0$ by the martingale equation for $\widetilde{d}_j$. Thus in this case we also have 
$\widetilde{d}_j(x) = \delta_j (x) \widetilde{d}_j (pref(x))$.

Note also that $\delta_j$ takes values in $\{0,1,2\}$.

We set 
 
\[  \widehat{d}_j(x)  = 
      \begin{cases} 
        \text{the constant value }  2^{-j} \   &\text{if} \  |x| < j \\ 
        \delta_j (x) \widehat{d}_j (pref(x)) \ &\text{if} \ |x| \ge j 
       \end{cases} 
\] 
 
From the martingale equation for $\widetilde{d}_j$ we have: 
\[ 2\widetilde{d}_j (x) = \widetilde{d}_j (x0) + \widetilde{d}_j (x1) 
   =[ \delta_j (x0) +  \delta_j (x1) ] \widetilde{d}_j (x), \] 
and so $\delta_j (x0) +  \delta_j (x1) = 2$ for all $x$. Thus for 
$|x|\ge j-1$ we have 
\[  \widehat{d}_j (x0) + \widehat{d}_j (x1) =  \delta_j (x) \widehat{d}_j (x) 
         + \delta_j (x) \widehat{d}_j (x) = 2 \widehat{d}_j(x).\] 
Similarly, it follows from the definition that for $|x|<j-1$ we have $\widehat{d}_j (x0) + \widehat{d}_j (x1)= 2 \widehat{d}_j(x)$. 
Thus $\widehat{d}_j$ satisfies the martingale equation.

It is easy to see that  $\widehat{d}_j $ 
 succeeds on $\widetilde{S}_j$.  To calculate  $\widehat{d}_j $ 
we need to compute $\widetilde{d}_j (y)$ and $\delta_j (y)$ on the 
prefixes $y$ of $x$  and this can be done  in time $n^{c_5}$ on 
inputs $x$ with $|x|=n \ge j$.

 We put these martingales together in the ``global'' martingale 
\begin{align} 
  \widehat{d} (x) &=  \sum_{j = 1}^{\infty} \widehat{d}_j (x) \\ 
  &= \sum_{j = 1}^{|x|}  \widehat{d}_j (x) +  \sum_{j = |x| + 1}^{\infty} 2^{-j} \\ 
  &= \sum_{j = 1}^{|x|} \widehat{d}_j (x) + 2^{-|x|} 
\end{align}

Then $\widehat{d}$ is a martingale which is calculable in 
polynomial time by a Turing machine with an oracle for $A$.  For 
each $j$ and $x$ with $|x| \ge j$ we have $\widehat{d}(x) > 
\widehat{d}_j (x)$, so since $ \widehat{d}_j (x)$ succeeds on 
$\widetilde{S}_j$ then $\widehat{d}$ succeeds on $\widetilde{S}_j$. 
This implies that $\widehat d$ succeeds on $\displaystyle 
\mathcal{UFP}^A=\cup_j \widetilde{S}_j=\cup_ i S_i$ and hence 
$\mathcal{UFP}^A$ has effective $PA$-measure zero, as claimed.

\end{proof} 
 
\begin{cor} We have: 
  \[\mathcal{NP}^A - \mathcal{UFP}^A \ne \emptyset.\] 
\end{cor} 
 
Thus partial complexity cannot compensate for nondeterminism in the 
presence of a random oracle and it is reasonable to suppose that some 
similar separation remains true without an oracle. For example, let 
$\mathcal{GP}$ be the class of languages which are generically 
decidable in polynomial time.  The assumption that $\mathcal{NP} - 
\mathcal{GP} \ne \emptyset$, would say that there are languages in 
$\mathcal{NP}$ which require nondeterminism on a nonnegligible set of 
inputs and is certainly a stronger hypothesis than just assuming that 
$\mathcal{NP} \ne \mathcal{P}$. It would be interesting to investigate 
the question of whether such a quantitative hypothesis yields stronger 
consequences.

 \end{document}